\documentclass[11pt]{article}

\usepackage{amsmath,amsthm,verbatim,amssymb,amsfonts,amscd, graphicx}
 \usepackage[matrix,arrow,ps]{xy}
\usepackage{graphics,pb-diagram}
\theoremstyle{plain}
\newtheorem{theorem}{Theorem}
\newtheorem{corollary}{Corollary}
\newtheorem{lemma}{Lemma}
\newtheorem{proposition}{Proposition}

\newtheorem*{theoremmain}{Theorem~\ref{thm:main}} 
\newtheorem*{theoremtriple}{Theorem~\ref{thm:triple}}

\theoremstyle{definition}
\newtheorem{definition}{Definition}

\newcommand{\Z}{\mathbb{Z}}
\newcommand{\R}{\mathbb{R}}

\newcommand{\zp}{\mathbb{Z}_p}
\newcommand{\A}{\widetilde{\mathcal{A}}}
\newcommand{\K}{\widetilde{K}}

\begin{document}

\title{Invariants for Legendrian knots in lens spaces}
\author{Joan E. Licata\\
\small\textit{Max-Planck-Institute f\"{u}r Mathematik/}\\
\small\textit{Stanford University}\\
\small\textit{jelicata@stanford.edu}}
\maketitle


\begin{abstract}
In this paper we define invariants for primitive Legendrian knots in lens spaces $L(p,q), q\neq 1$.  The main invariant is a differential graded algebra $(\mathcal{A}, \partial)$ which is computed from a labeled Lagrangian projection of the pair $(L(p,q), K)$.  This invariant is formally similar to a DGA defined by Sabloff which is an invariant for Legendrian knots in smooth $S^1$-bundles over Riemann surfaces.   The second invariant defined for $K\subset L(p,q)$ takes the form of a DGA enhanced with a free cyclic group action and can be computed from the $p$-fold cover of the pair $(L(p,q), K)$.

\end{abstract}

\section{Introduction}
Endowing a three-manifold with a contact structure refines the associated knot theory by introducing new notions of equivalence among knots, and these in turn require invariants sensitive to the added geometry.  In addition to more classical numerical invariants, invariants taking the form of differential graded algebras (DGAs)  have seen success in distinguishing Legendrian non-isotopic knots in a variety of contact manifolds.  

The first DGA invariants were developed for Legendrian knots in the standard contact $\mathbb{R}^3$.  Chekanov constructed a combinatorial invariant, and an equivalent invariant was introduced independently in a geometric context by Eliahsberg \cite{C}, \cite{E}.  In the former case, the algebra is generated by the crossings in a Lagrangian projection, and the boundary map counts immersed discs in the diagram.  In Eliashberg's relative contact homology, the algebra is generated by Reeb chords and the differential counts rigid $\mathcal{J}$-holomorphic curves in the symplectization of $\R^3$.   The two constructions were shown to produce the same DGA in \cite{ENS}.  In \cite{S}, Sabloff adapted further work of Eliashberg, Givental, and Hofer to construct a combinatorial DGA for Legendrian knots in a class of contact manifolds characterized by their distinctive Reeb dynamics \cite{EGH}, \cite{S}.  His algebra is again generated by Reeb chords, but he introduces additional technical machinery in order to handle periodic Reeb orbits.  Sabloff's invariant is defined for smooth $S^1$ bundles over Riemann surfaces, a class of manifolds which includes $S^3$ and $L(p,1)$, but does not admit other lens spaces.

In this paper we develop an invariant for Legendrian knots in the lens spaces $L(p,q)$ for $q \neq 1$ with the unique universally tight contact structure.  For primitive $K \subset L(p,q)$, we define a \textit{labeled diagram} to be the Lagrangian projection of the pair $(L(p,q), K)$ to $(S^2, \Gamma)$, together with some ancillary decoration which uniquely identifies the Legendrian knot.  Numbering the crossings of $\Gamma$ from one to $n$, we consider the tensor algebra on $2n$ generators:
\[\mathcal{A}=T(a_1, b_1, ...a_n, b_n).
\]
We equip this algebra with a differential $\partial:\mathcal{A}\rightarrow \mathcal{A}$ counting certain immersed discs in $(S^2, \Gamma)$.  The algebra is graded by a cyclic group, and the boundary map is graded with degree $-1$.  The pair $(\mathcal{A}, \partial)$ is a semi-free DGA, and the natural equivalence on such pairs is that of stable tame isomorphism type.  

Our main theorem is the following:
\begin{theoremmain} Up to equivalence, the semi-free DGA $(\mathcal{A}, \partial)$  is an invariant of the Legendrian type of $K\subset L(p,q)$. 
\end{theoremmain}

The proof of Theorem~\ref{thm:main} applies Sabloff's invariant to a freely periodic knot $\K \subset S^3$ which is a $p$-to-one cover of $K\subset L(p,q)$.  The Legendrian type of $\K$ is an invariant of the Legendrian type of $K$, so Sabloff's invariant  for $\K$ is therefore also an invariant of $K$  (Proposition~\ref{prop:easyinvt}).  In order to prove Theorem~\ref{thm:main}, we endow Sabloff's DGA with additional structure related to the covering transformations. 

Given $K$ in $L(p,q)$ with $q\neq1$, let $(\widetilde{\mathcal{A}}, \widetilde{\partial})$ denote Sabloff's low-energy DGA for the knot $\K$ in $S^3$.
The algebra $(\widetilde{\mathcal{A}}, \widetilde{\partial})$ may be enhanced with a cyclic group action $\gamma:\zp \times \widetilde{\mathcal{A}} \rightarrow \widetilde{\mathcal{A}}$ which commutes with the boundary map.  We define a notion of equivariant equivalence on DGAs with such actions in Section~\ref{sect:triplepf}, and we associate to $K$ the equivariant DGA $(\widetilde{\mathcal{A}}, \gamma, \widetilde{\partial})$.  Our second main theorem asserts that this is also an invariant of the Legendrian knot in the lens space.

\begin{theoremtriple} The equivalence class of the equivariant DGA $(\widetilde{\mathcal{A}}, \gamma, \widetilde{\partial})$ is an invariant of the Legendrian type of $K$.
\end{theoremtriple}

The major technical work of the paper lies in proving Theorem~\ref{thm:triple}, and this occupies Section~\ref{sect:triplepf}.  The proof of Theorem~\ref{thm:main} identifies $(\mathcal{A}, \partial)$ with a distinguished $\zp$-equivariant subalgebra of $(\widetilde{\mathcal{A}}, \widetilde{\partial})$ and follows as a consequence of Theorem~\ref{thm:triple}.   The final section contains examples computed for knots in $L(3,2)$ and $L(5,2)$.

Finally, we note that although the arguments in this paper are developed for primitive knots in $L(p,q)$, they in fact construct invariants for any Legendrian knot in a lens space  which is covered by a Legendrian knot in some $L(p,1)$.  In this adaptation, $L(p,1)$ replaces $S^3$ as the contact manifold where Sabloff's invariant is defined.

I would like to thank Josh Sabloff for helpful correspondence in the course of writing this paper.   A portion of this work was conducted while visiting the Max Planck Institute for Mathematics in Bonn, Germany, and their hospitality and support are much appreciated.

\section{Background}\label{background}
This section contains a brief summary of the basic definitions from contact geometry and their realizations in three examples: the standard contact $\R^3$, $S^3$, and $L(p,q)$.  A more thorough introduction to the topic is provided in \cite{Et} or \cite{Ge}.

\subsection{Basic definitions}
A \textit{contact structure} $\xi$ on a three-manifold $M$ is an everywhere non-integrable $2$-plane field.  A non-degenerate one-form $\alpha$ defines a contact structure by $\xi_{\alpha}= \ker \alpha$ at each point of $M$.  Two contact manifolds $(M_1, \xi_1)$ and $(M_2, \xi_2)$ are \textit{contactomorphic} if there is a diffeomorphism between the manifolds which takes contact planes to contact planes.

\begin{definition} 
Given a contact form $\alpha$, the \textit{Reeb vector field}  is the unique vector field $X$ which satisfies 
\begin{align}
&\alpha(X)=1\notag \\
& d\alpha(X, \cdot)=0. \notag
\end{align}
Integral curves of $X$ are known as \textit{Reeb orbits}, and they inherit an orientation from $X$.
\end{definition}

\begin{definition} 
A knot $K$ in $(M, \xi)$ is \textit{Legendrian} if its tangent lies in the contact plane at each point.
\end{definition}

Two Legendrian knots are equivalent if they are isotopic through Legendrian knots.  In general, two knots which are topologically equivalent may not be Legendrian equivalent; any topological isotopy class of knots will be represented by countably many Legendrian isotopy classes.

\begin{definition} 
The \textit{Lagrangian projection} of a contact manifold $(M, \xi_{\alpha})$ is the quotient space of $M$ which collapses each Reeb orbit of $\alpha$ to a point.  If $K$ is a Legendrian knot in a contact manifold, the \textit{Lagrangian projection} of $K$ is the image of the knot under Lagrangian projection of the manifold.
\end{definition}

If $K$ is a Legendrian knot in $(M, \xi)$, a \textit{Reeb chord} is a segment of a Reeb orbit with both endpoints on $K$.  In the Lagrangian projection, a Reeb chord with distinct endpoints will map to a crossing in the knot projection.

\subsection{First example: $\R^3$} The standard contact structure $\xi_{std}$ on $\R^3$ is induced by the contact form \[
\alpha_{std}=dz-ydx.
\]   The Reeb vector field on $(\R^3, \xi_{std})$ has trivial $dx$ and $dy$ coordinates at every point, so the Reeb orbits are vertical lines.  Thus, the Lagrangian projection is simply projection to the $xy$-plane.

\subsection{Second example: $S^3$}\label{sect:s3}
$S^3$ sits inside $\R^4$ as the unit sphere:
\[S^3=\{ (r_1, \theta_1, r_2, \theta_2) | r_1^2+r_2^2=1\}.
\]
The torus $r_1=\frac{1}{\sqrt{2}}=r_2$ separates $S^3$ into two solid tori, and it will be convenient to treat this torus as a Heegaard surface.  The curves $r_1=0$ and $r_2=0$ are the core curves of the Heegaard tori, and the complement of the cores is  foliated by tori of fixed $r_i$.  

The standard tight contact structure on $S^3$ is 
\[
\alpha_0=\frac{1}{2}(r_1^2d\theta_1+r_2^2 d\theta_2).
\]
The punctured manifold $(S^3 - \{p\}, \xi_0)$ is contactomorphic to $(\R^3, \xi_{std})$, but the Reeb dynamics are quite different.  In particular, the Reeb orbits of $\alpha_0$ are $(1,1)$ curves on each torus of fixed $r_i$.  This foliation of $S^3$ by circles gives the Hopf fibration of $S^3$, and Lagrangian projection in $(S^3, \xi_0)$ is projection to the $S^2$ base space of this fibration. Note that the core curves are each Reeb orbits, and their images under Lagrangian projection are the poles of the two-sphere.  The contact form $\alpha$ also induces a \textit{curvature} form $\Omega$ on the
 $S^2$ base space;  for the standard contact structure, this is just the Euler class of the bundle, where $S^3$ is viewed as the unit sphere in $\mathbb{R}^4$. \cite{Ge}.    

\subsection{Third example: Lens spaces}\label{sect:lens}
Define $F_{p,q}:S^3 \rightarrow S^3$ by
\begin{equation}\label{eq:fp}
F_{p,q}(r_1, \theta_1, r_2, \theta_2)=(r_1, \theta_1+ \frac{2 \pi}{p}, r_2,\theta_2 + \frac{2 q\pi}{p}).
\end{equation}

The map $F_{p,q}$ generates a cyclic group of order $p$, and the quotient of $S^3$ by the action of this group is the lens space $L(p,q)$.  Thus $\pi:S^3\rightarrow L(p,q)$ is a $p$-to-one covering map.  Since $F_{p,q}$ preserves the contact structure on $S^3$, $\pi$ induces a contact structure on $L(p,q)$ \cite{BG}.  The Reeb orbits of $(L(p,q), \xi_{p,q})$ again foliate the manifold by circles, and the Lagrangian projection of $(L(p,q), \xi_{p,q})$ is a two-sphere.  As an $S^1$ bundle over $S^2$,  $L(p,q)$  is smooth if and only if $q=1$. 

\begin{definition}\label{def:fp} 
A knot $\K$ in $S^3$ is \textit{freely periodic} if it is preserved by a free periodic automorphism of $S^3$.
\end{definition}
 
  The map in Equation \ref{eq:fp} has order $p$, so if $\K$ is freely periodic with respect to $F_{p,q}$, then $\pi(\K)$ is a knot in  $L(p,q)$.  Conversely, any $K$ in $L(p,q)$ which is primitive in $H_1(L(p,q))$ has a freely periodic lift $\K\subset S^3$.  (Knots which are not primitive will lift to links in $S^3$.)  This definition makes sense in both the topological and contact categories; with respect to the contact structures defined above, $K$ is Legendrian if and only if $\K$ is Legendrian.  An explicit construction of a freely periodic lift is described in Section 6.2 of \cite{GRS}, and we refer the reader to \cite{HLN} or \cite{R2} for a fuller treatment of freely periodic knots. 

Throughout the paper, each topological manifold will be equipped with the contact structure associated to it in this section; we will write only $S^3$ and $L(p,q)$ for the contact manifolds $(S^3, \xi_{std})$, and $(L(p,q), \xi_{p,q})$.   Furthermore, tildes will be used to distinguish objects in $S^3$ from their counterparts in $L(p,q)$; thus $\widetilde{\Gamma}$ will denote the Lagrangian projection of a knot $\K$ in $S^3$, whereas the Lagrangian projection of $K\subset L(p,q)$ will be denoted by $\Gamma$.

\section{Differential graded algebra invariants for Legendrian knots}\label{DGA}
In this section we introduce Sabloff's DGA invariant for Legendrian knots in smooth $S^1$ bundles over Riemann surfaces.  We begin by defining differential graded algebras and the relevant notion of equivalence among them.

\subsection{Equivalence of semi-free DGAs}\label{sect:eq}

Let $V=\text{Span}_{\mathbb{Z}_2} \{x_1, x_2, ...x_n\}$.  Define
\[
\mathcal{A}= T(x_1, x_2, ...x_n)=\bigoplus_{n=0}^{\infty} V^{\otimes n}
\]
to be the tensor algebra on the elements $\{x_1, x_2, ...x_n\}$.  If $V$ is graded by a cyclic group $G$ so that the $x_i$ are homogeneous, this induces a cyclic grading on $\mathcal{A}$ via the rule $|x_ix_j|=|x_i|+|x_j|$.  When $\partial:\mathcal{A}\rightarrow \mathcal{A}$ is a degree $-1$ map satisfying  $\partial^2=0$ and the Leibnitz rule $\partial (ab)=(\partial a)b+a(\partial b)$, then the pair $(\mathcal{A}, \partial)$ is a \textit{semi-free differential graded algebra} (DGA).  The modifier ``semi-free" emphasizes that we keep track of the preferred generators $\{x_i\}_{i=1}^n$, which will be important in defining DGA equivalence.

An \textit{elementary automorphism} of $\mathcal{A}$ is a map $g^i: \mathcal{A}\rightarrow \mathcal{A}$ such that
\[
g^i(x_j)= \begin{cases}  x_i+ v_i, \text{ for } v_i \in T(a_1,....\hat{x}_i ...b_n) &  \text{if }j=i  \\ 
x_j & \text{if }  j\neq i. \end{cases} 
\]
 When $v_i$ is homogeneous in the same grading as $x_i$, we say that $g^i$ is a \textit{graded elementary automorphism}.  A \textit{graded tame automorphism} is a composition of graded elementary automorphisms.  

Given a DGA $(\mathcal{A}, \partial)=(T(x_1, ...x_n), \partial)$, let  $\mathcal{E}=T(e_1, e_2)$ be a DGA which is graded by the same cyclic group and satisfies  $\partial_{\mathcal{E}}e_1=e_2$ and $\partial_{\mathcal{E}}e_2=0$.  A \textit{stabilization} of $\mathcal{A}$ is the differential graded algebra $(T(x_1, ...x_n, e_1, e_2),  \partial \coprod \partial_{\mathcal{E}})$.  

\begin{definition}\label{def:equiv}
Two semi-free differential graded algebras $(\mathcal{A}_1, \partial_1)$ and $(\mathcal{A}_2, \partial_2)$ are \textit{equivalent} if some stabilization of $(\mathcal{A}_1, \partial_1)$ is graded tame isomorphic to some stabilization of $(\mathcal{A}_2, \partial_2)$.
\end{definition}

We will have reason to consider DGAs equipped with an action of a cyclic group $\zp$, so we extend the notion of equivalence to one respecting the group action.  The cyclic group  $\zp$ should not be confused with the cyclic group $G$ which grades the algebra. 

\begin{definition} An \textit{equivariant DGA} $(\mathcal{A}, \gamma, \partial)$ is a semi-free DGA $(\mathcal{A}, \partial)$ together with an automorphism $\gamma: \mathcal{A} \rightarrow \mathcal{A}$ of order $p$  such that $\partial \circ \gamma=\gamma \circ \partial$ and $|\gamma x|=|x|$.
\end{definition}

\begin{definition}  Suppose that  $(\mathcal{A}, \gamma, \partial)$ is an equivariant DGA, where $\mathcal{A}=T(x_1, ...x_n)$.  A \textit{free} $\zp$ \textit{stabilization} of $(\mathcal{A}, \gamma, \partial)$ is the equivariant DGA $(\mathcal{A}\coprod \mathcal{E}^p, \gamma, \partial \coprod \partial_{\mathcal{E}^p})$, where 
\begin{itemize}
\item $\mathcal{A}\coprod \mathcal{E}^p = T(x_1, ...x_n, e_{1,1}, e_{1,2}, ...e_{1,p}, e_{2,1}, ...e_{2,p})$;
\item $\gamma(e_{j,i})=e_{j, i+1}$;
\item $\partial_{\mathcal{E}^p} e_{1,i}=e_{2,i}$;
\item $\partial_{\mathcal{E}^p} e_{2,i}=0$.
\end{itemize}
\end{definition}

\begin{definition} A $\zp$ \textit{elementary isomorphism} is a $\gamma$-equivariant map $f: (\mathcal{A}, \gamma, \partial)\rightarrow (\mathcal{A}, \gamma, \partial)$ which can be written as $f=g^1\circ g^2\circ...\circ g^p$, where each $g^i$ is an elementary isomorphism. If $f$ is graded, we say it is a \textit{graded} $\zp$ \textit{elementary isomorphism}.  A composition of $\zp$ elementary isomorphisms is a $\zp$ \textit{tame isomorphism}.
\end{definition}

\begin{definition} Two equivariant DGAs $(\mathcal{A}_1, \gamma_1, \partial_1)$ and $(\mathcal{A}_2, \gamma_2, \partial_2)$ are $\zp$ \textit{equivalent} if they have free $\zp$ stabilizations which are graded $\zp$ tamely isomorphic.  \end{definition}

\subsection{Sabloff's DGA for knots in $S^1$ bundles over Riemann surfaces}\label{sect:sab}
In \cite{S}, Sabloff considers contact manifolds whose Reeb orbits are the fibers of a smooth $S^1$ bundle over a Riemann surface.  For a Legendrian knot $\K$ in such a manifold, he defines  an algebra generated by the Reeb chords with both endpoints on $\K$.  Since each Reeb orbit is periodic, there are infinitely many such chords, and he also defines a finitely-generated \textit{low-energy algebra} generated by chords which are strictly shorter than the fiber.  The low-energy algebra sits inside the full invariant as a subalgebra, but the equivalence type of the low-energy algebra is also an invariant of $\K$.  The following section introduces the the low-energy algebra $(\A, \K)$ for Legendrian knots in $S^3$, and we refer the reader to \cite{S} for a description of the full invariant.
  
\subsubsection{Labeled Lagrangian diagram}

\begin{figure}[h]
\begin{center}
\scalebox{.7}{\includegraphics{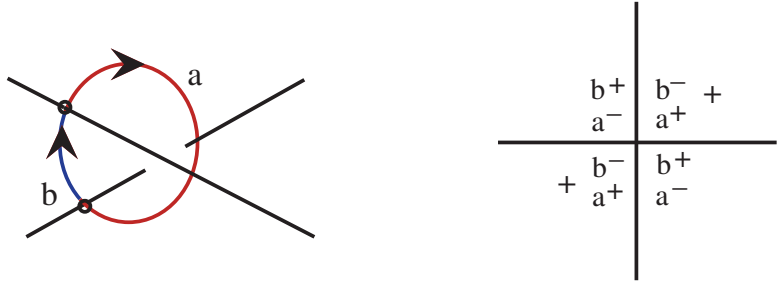}}
\end{center}
\caption{Labels at a crossing on $\Gamma$. The ``+" in the right-hand diagram indicates that $a$ is the preferred chord.}\label{fig:crosslabels}
\end{figure}

Let $\K\subset S^3$ be a Legendrian knot, and denote the Lagrangian projection of $\K$ to $S^2$ by $\widetilde{\Gamma}$.  Number the crossings of the diagram from $1$ to $n$, and associate two generators  $a_i$ and $b_i$ to the $i^{th}$ crossing.  These correspond to the complementary short chords in the fiber which intersects the crossing strands of $\K$.  Because the Reeb orbit is oriented, each chord identifies the crossing strands locally as ``sink" and ``source".  Select a preferred chord and indicate this choice with a plus sign in the two (opposite) quadrants where traveling sink-to-source orients the quadrant positively.   Furthermore, assign each quadrant either $a_i^+$ and $b_i^-$ or $a_i^-$ and $b_i^+$ as indicated in Figure~\ref{fig:crosslabels}.  Note that signs are used in two distinct ways; a ``positive quadrant" will always mean one marked with a ``+" to denote the preferred chord, and each generator $x$ labels every quadrant of the associated crossing as either $x^+$ or $x^-$.
  
\begin{definition} If $\widetilde{\Gamma}$ is a labeled diagram with $n$ crossings, define $\A(\widetilde{\Gamma})$ to be the tensor algebra generated by the associated Reeb chords:
\[
\A(\widetilde{\Gamma})=T(a_1, b_1, ...a_n, b_n).
\]
\end{definition}

The following definitions will prove useful in defining the defect and the boundary map:
 \begin{definition}
If $x_i$ is a generator of $\A(\widetilde{\Gamma})$, let $\mathit{l}(x_i)$ denote the length of the associated chord in $S^3$, where the length of an $S^1$ fiber is normalized to $1$.  
\end{definition}

We extend this to a length function $\mathit{l}'$ on words written in the signed generators $a_i^{\pm}$ and $b_i^{\pm}$.  Let $\epsilon(a_i^+)=\epsilon(b_i^+)=1$ and $\epsilon(a_i^-)=\epsilon(b_i^-)=-1$.  If $w$ is a word in the signed generators $x_i^{\pm}$, define
\[
\mathit{l}'(w)=\sum_{x_i^{\pm} \in w} \epsilon(x_i^{\pm}) \mathit{l}(x_i).
\]

\begin{definition}\label{def:adm}
Let $(\Sigma, \partial \Sigma)$ be a disc with $m$ marked points on the boundary.
An \textit{admissible disc} is a map $f:(\Sigma, \partial \Sigma) \rightarrow (S^2, \widetilde{\Gamma})$  which satisfies the following:
\begin{enumerate}
\item each marked point maps to a crossing of $\widetilde{\Gamma}$;
\item $f$ is an immersion on the interior of $\Sigma$;
\item $f$ extends smoothly to $\partial \Sigma$ away from the marked points;
\item $f(\partial \Sigma)$ has a corner at each marked point, and $f(\Sigma)$ fills one quadrant there.
\end{enumerate}
\end{definition}

Two admissible discs $f$ and $g$ are \textit{equivalent} if there is a smooth automorphism $\phi: \Sigma \rightarrow \Sigma$ such that $f=g\circ \phi$.

\begin{figure}[h]
\begin{center}
\scalebox{.5}{\includegraphics{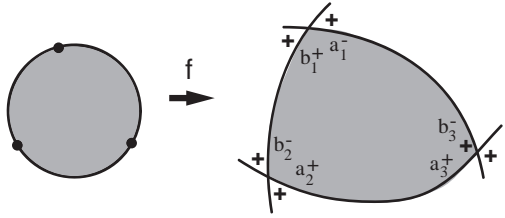}}
\end{center}
\caption{An admissible disc with $w(R)=a_1^-b_2^-a_3^+$.  This disc could represent three different boundary words: $w(f, b_1)=b_2 b_3$; $w(f,a_2)=b_3 a_1$; or $w(f, a_3)=a_1 b_2$.}\label{fig:word}
\end{figure}

Let $R$ be a component of  $S^2-\widetilde{\Gamma}$.  To each corner of $R$, one may associate the signed generator corresponding to the preferred chord of the crossing, where the sign is dictated by the quadrant filled by $R$.  Traveling counterclockwise around $\partial R$ and reading off these labels defines a cyclic word $w(R)$. (See Figure~\ref{fig:word} for an example.)

 \begin{definition}\label{def:def}
   Let $f$ be an admissible disc whose image is $R$.  The \textit{defect} of $R$ is given by: \begin{equation}\label{eq:def}
 n(f)=\frac{1}{2\pi}\int_{\Sigma} f^*\Omega + \mathit{l}'(w(R)).
 \end{equation}
 \end{definition}
 
Geometrically, the defect encodes the interaction between the knot and the fiber structure.  Without this decoration, $\widetilde{\Gamma}$ does not specify even the topological type of the knot, as displacement in the Reeb direction is obscured by the projection.  The curve $\partial R$ lifts to $S^3$ as a simple closed curve composed of alternating Legendrian and Reeb segments, and the defect measures the winding number of this lifted curve around the fiber with respect to an appropriate trivialization. Together with the signs at each crossing, the defects of components of $S^2-\widetilde{\Gamma}$ determine the Legendrian type of the knot.  

 \begin{figure}[h]
\begin{center}
\scalebox{.4}{\includegraphics{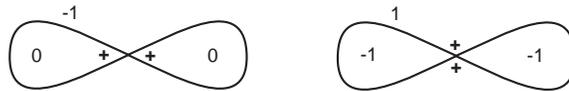}}
\end{center}
\caption{Two labeled diagrams for the same Legendrian unknot in $S^3$.}\label{fig:unknots}
\end{figure}

The defect extends additively to unions of regions counted with multiplicity, so  Equation~\ref{eq:def} holds for any admissible disc $f$. Since $S^2$ is simply connected, one may also define the defect of the knot $n(\K)$ to  be the defect of any contracting disc bounded by the projection of $\K$.   

\subsubsection{Gradings}

A \textit{capping path} for a generator $x_i$ is a path along $\K$ in $\widetilde{\Gamma}$ which begins and ends adjacent to the same $x_i^+$ quadrant.  For each crossing, one of $a_i$ or $b_i$ will have two capping paths, and the other will have none.  The \textit{rotation number} of a capping path for $x_i$ is the number of counterclockwise rotations performed by the tangent vector, computed as a winding number in a trivialization over a contracting disc in $S^2$.  Taking the edges at a crossing to be orthogonal, this value lies in  ${\Z}-\frac{1}{4}$, and we denote it by $r(x_i)$.  

Suppose that $f:\Sigma \rightarrow S^2$ is an admissible disc such that $f(\partial \Sigma )$ is a capping path for $x_i$ .  Then the grading of $x_i$ is given by the following:
\begin{equation}\label{eq:grad1}
|x_i|= 2 r(x_i) -\frac{1}{2}+4n(f).
\end{equation}
If $y_i$ is the other generator at the same crossing,
\begin{equation}\label{eq:grad2}
|y_i|= 3-|x_i|.
\end{equation}

These gradings are well-defined modulo $2r(x_i) + 4n(\K)$.

\subsubsection{The algebra $(\A(\widetilde{\Gamma}), \widetilde{\partial})$}

\begin{definition}\label{def:word}
  Let $f:\Sigma\rightarrow S^2$ be an admissible disc with one corner filling a quadrant labeled $x_i^+$.  The \textit{boundary word} $w(f,x_i)$ is the concatenation of the $y_j^-$ generators associated to the other quadrants filled by $f(\Sigma)$, read counterclockwise around $\partial \Sigma$.  
\[w(f,x_i)=y_2 y_3...y_{m}.\]
\end{definition}

See Figure~\ref{fig:word} for an example.

\begin{definition}
If $f:\Sigma\rightarrow S^2$ is an admissible disc, the \textit{$x_i$ defect}  $\tilde{n}_{x_i}(f)$ is given by
\[
\tilde{n}_{x_i}(f)=\frac{1}{2\pi}\int_{\Sigma} f^*\Omega + \mathit{l}(x_i)-\sum_{j=2}^{m} \mathit{l}(y_j).
\]
\end{definition}

Note that the $x_i$ defect of an admissible disc may differ from the defect of its image in the diagram, as the two are computed by associating (possibly) different words to the same disc. To compute $\tilde{n}_{x_i}(f)$, add one to $n(f)$ if $x_i^+$ occupies a non-positive quadrant, and subtract one from $n(f)$ for each $y_j^-$ in $w(f, x_i)$ which occupies a positive quadrant.  Thus both types of defects may be computed from the labeled diagram without further data regarding the lengths of chords.

\begin{definition}\label{def:bndry}
The differential $\widetilde{\partial}:\A(\widetilde{\Gamma})\rightarrow \A(\widetilde{\Gamma})$ is defined by on the generator $x_i$  by
\[
\widetilde{\partial} x_i=\sum_{f: \tilde{n}_{x_i}(f)=0} w(f,x_i), 
\]
and $\widetilde{\partial}$ extends to other elements in $\A(\widetilde{\Gamma})$ via the Leibnitz rule $\widetilde{\partial}(ab)=(\widetilde{\partial}a)b+a(\widetilde{\partial}b)$.
\end{definition}

\begin{theorem}[\cite{S} Proposition 3.8, Theorem 3.11, Corollary 3.16]\label{sabinvt}
The boundary map in Definition~\ref{def:bndry} satisfies $\widetilde{\partial}^2=0$, and the stable tame isomorphism type of $(\A(\widetilde{\Gamma}), \widetilde{\partial})$ is an invariant of the Legendrian knot type of $\K$ in $S^3$.
\end{theorem}

\section{Invariants for Legendrian knots in lens spaces}\label{sect:invt}

As noted above, Sabloff's invariant is defined for contact manifolds which are smooth $S^1$ bundles, a class which excludes the lens spaces $L(p,q)$ for $q\neq1$.  Although they do not induce smooth bundles, the Reeb orbits of these lens spaces nevertheless define an $S^1$ bundle structure, and this similarity is strong enough to permit a DGA invariant $(\mathcal{A},\partial)$ computable from the Lagrangian projection to $S^2$.  The invariant is formally similar to Sabloff's invariant $(\A, \widetilde{\Gamma})$ for knots in $S^3$, and in fact, the proof of invariance exploits the covering relationship between these manifolds.  
 
Except if otherwise indicated, in the remainder of the paper every lens space $L(p,q)$ is assumed to have $q\neq1$.
 
 \subsection{The DGA $(A, \partial)$}
 
 Let $K$ be a knot in $L(p,q)$ which generates $H_1(L(p,q))$.  Following Rasmussen, we call such knots \textit{primitive} \cite{R2}.  If $K$ is a primitive Legendrian knot, we begin by defining a labeled Lagrangian diagram.  At the $i^{th}$ crossing of $\Gamma$, mark each quadrant with $a_i^+$ and $b_i^-$ or with $a_i^-$ and $b_i^+$ as in Figure~\ref{fig:crosslabels}.  At each crossing, indicate a preferred choice of chord by decorating a pair of opposite quadrants with plus signs.  

Recall that to $K$, we may associate its freely perioidic lift $\K\subset S^3$.  The $p$-fold covering map $\pi:(S^3, \K)\rightarrow (L(p,q), K)$ descends to a $p$-to-one branched cover of Lagrangian projections $\pi_*:(S^2, \widetilde{\Gamma})\rightarrow (S^2, \Gamma)$, where the branch points are the images of the core curves $r_i=0$ for  $i=1,2$.  Thus a choice of preferred chords in $\Gamma$ lifts to a choice of preferred chords in $\widetilde{\Gamma}$.  Let $R$ be a region in $(S^2 - \Gamma)$, and let $f:\Sigma\rightarrow S^2$ be an admissible disc whose image is $\pi_*^{-1}(R)$.  Define the \textit{defect} of $R$ to be $n(R)=\frac{1}{p}n(f)$.

A \textit{labeled diagram} for $K$ is a generic Lagrangian projection $\Gamma$ decorated with preferred chords and defects which are compatible with a labeled diagram for $\K$ as described above.

\begin{definition} Let $\Gamma$ be a labeled diagram for a Legendrian knot $K\subset L(p,q)$.  If $\Gamma$ has $n$ crossings, define
\[
\mathcal{A}(\Gamma)=T(a_1, b_1, ...a_n, b_n).
\]
\end{definition}

If $x_i$ is a generator of $\mathcal{A}(\Gamma)$, choose a lift $x_i^*\in \pi_*^{-1}(x_i)$ and define the grading of  $x_i$ by  $|x_i|=|x_i^*|$.  This value is independent of the choice of lift, and $\mathcal{A}(\Gamma)$ is graded by the same cyclic group as $\A(\widetilde{\Gamma})$.  

\textbf{Remark} The grading can also be defined intrinsically.  Given a labeled diagram, consider a capping path for $x_i$ which has winding number $p$ with respect to the poles.  With only slight modification,  the formulae in Equations~\ref{eq:grad1} and \ref{eq:grad2} can be used to compute the grading directly from $\Gamma$.  

Definitions~\ref{def:adm} and \ref{def:word} may be applied verbatim in the context of labeled diagrams for knots in $L(p,q)$.

\begin{definition}\label{def:bndry2}
The differential  $\partial:\mathcal{A}(\Gamma)\rightarrow \mathcal{A}(\Gamma)$ is defined on generators by
\[
\partial x_i=\sum_{f: \tilde{n}_{x_i}(f)=0} w(f,x_i),
\]
where the sum is over admissible discs which satisfy the additional condition that $f(\partial \Sigma)$ has winding number $p$ with respect to the poles of $S^2$.  Extend $\partial$ to other elements in $\mathcal{A}(\Gamma)$ via the Leibnitz rule.
\end{definition}

\begin{theorem}\label{thm:main}
Up to equivalence as a semi-free DGA, $(\mathcal{A}(\Gamma), \partial)$ is an invariant of the Legendrian knot type of $K\subset L(p,q)$.
\end{theorem}
 
In order to prove Theorem~\ref{thm:main}, we will study the relationship between $(\mathcal{A}(\Gamma), \partial)$ and $(\A(\widetilde{\Gamma}), \widetilde{\partial})$.  
 
It is clear that the Legendrian type of the freely periodic lift $\K \subset S^3$ is an invariant of the Legendrian type of $K$, so Sabloff's construction has the following easy consequence: 

\begin{proposition}\label{prop:easyinvt}  The stable tame isomorphism type of $(\A(\widetilde{\Gamma}), \widetilde{\partial})$ is an invariant of the Legendrian isotopy class of $K$.  
\end{proposition}

However, a stronger notion of equivalence yields a more interesting invariant, and the next section shows that we may associate an equivariant DGA to the freely-periodic lift of $K$.  

\subsection{The $\zp$ action on $(\mathcal{A}(\Gamma), \partial)$}\label{sect:zpact}

Theorem 3.15 of \cite{S} states that the equivalence type of $(\A(\widetilde{\Gamma}), \widetilde{\partial})$ is independent of the choice of preferred chords, but we will restrict attention to diagrams  where the signs at each crossing are preserved by $\frac{2 \pi}{p}$ rotation of $(S^2, \widetilde{\Gamma})$.  

 \begin{lemma} If $K$ is a Legendrian knot in $L(p,q)$, then there is a natural automorphism $\gamma: \A(\widetilde{\Gamma}) \rightarrow \A(\widetilde{\Gamma})$ with order $p$ such that  $\widetilde{\partial}\circ \gamma=\gamma \circ \widetilde{\partial}$, and $|\gamma x|=|x|$.
 \end{lemma}
 
\begin{proof} 
Fix a representative of the isotopy class of $K$ and lift this to the freely periodic knot $\K$.  The Lagrangian projection of $(S^3,\K)$ is invariant under $\frac{2 \pi}{p}$ rotation about the axis through the points representing the fibers $r_i=0$ for $i=1,2$. If the projection of $\K$ is not generic, any local perturbation of $K$ will lift to $p$ local perturbations of $\K$, maintaining the contact covering relationship between $(S^3, \K)$ and $(L(p,q),K)$ while removing  singularities in the projection.  In particular, $K$ (or equivalently, $\K$) may be assumed disjoint from the cores of the Heegaard tori.  

 If the Reeb chord $x$ is a generator of $\mathcal{A}$, then $\pi^{-1}(x)$ is a free $\zp$ orbit of generators of $\A(\widetilde{\Gamma})$.  This relationship descends to the Lagrangian diagrams, via the $p$-fold branched covering map $\pi_*:(S^2, \widetilde{\Gamma})\rightarrow(S^2, \Gamma)$.   Since capping paths for crossings in a single orbit are permuted by the cyclic action, each member of the orbit has the same grading.  Similarly, any disc which represents a term in the boundary is part of an orbit of $p$ discs.   This proves that the $\zp$ action on $\A(\widetilde{\Gamma})$ commutes with the differential. 
\end{proof}

\textbf{Remark} Recall that Sabloff's invariant is defined for knots in lens spaces $L(p,1)$.  The above discussion highlights another sense in which this case is exceptional.  When $q=1$, the map $\pi_*:(S^2, \widetilde{\Gamma})\rightarrow(S^2, \Gamma)$ induced on Lagrangian projections is one-to-one.  In this case, for any point $x$ on $\widetilde{\Gamma}$, the preimage $F^{-1}_{p,1}(x)$  consists of $p$ points on $\K$.    

\begin{theorem}\label{thm:triple}
 If $K$ is a Legendrian knot in $L(p,q)$, then the equivariant DGA $(\A(\widetilde{\Gamma}), \gamma, \widetilde{\partial})$ is an invariant of $K$, up to $\zp$ equivalence.
\end{theorem}

The proof of Theorem~\ref{thm:triple} appears later, but we note a corollary of the statement here:

\begin{corollary} Let $\mathcal{A}^{\gamma}(K)$ be the subalgebra of $(\A(\widetilde{\Gamma}, \gamma, \widetilde{\partial})$  fixed by the $\zp$ action:
\[
\mathcal{A}^{\gamma}(K)=\{ a \in \A(\widetilde{\Gamma}) | \gamma a=a \}.
\]
The subalgebra $\mathcal{A}^{\gamma}(K)$ is a subcomplex and the homology of $\mathcal{A}^{\gamma}(K)$ is an invariant of the Legendrian type of $K$.
\end{corollary}

\begin{proof} 
The statement that $\mathcal{A}^{\gamma}(K)$ is a subcomplex follows from the fact that $\gamma$ commutes with the differential.  A $\zp$-equivariant isomorphism between equivariant DGAs induces an isomorphism on their homologies; the proof is similar to that in $\cite{C}$.  The statement then follows  from Theorem~\ref{thm:triple}, 
\end{proof}

\subsection{Proof of Theorem~\ref{thm:main}}
In this section we will show how Theorem~\ref{thm:main} follows from Theorem~\ref{thm:triple}, postponing the proof of Theorem~\ref{thm:triple} until Section~\ref{sect:triplepf}.

\begin{definition}
Setting $\phi: \A(\widetilde{\Gamma}) \rightarrow \A(\widetilde{\Gamma})$ to be the algebra homomorphism defined on the generators of $\A(\Gamma)$ by
\[
\phi(x)=\bar{x}=\sum_{i=1}^{p} \gamma^i x, 
\]
define  $\bar{\mathcal{A}}(\widetilde{\Gamma})$ to be the image of $\A(\widetilde{\Gamma})$ under $\phi$. 
\end{definition}

$\bar{\mathcal{A}}(\widetilde{\Gamma})$ is a $\zp$-equivariant subalgebra of both $\A$ and $\mathcal{A}^{\gamma}(K)$, and for every Legendrian knot, the following containments are proper:
\[
\bar{\mathcal{A}}(\widetilde{\Gamma}) \subset \ \mathcal{A}^{\gamma}(K) \subset \A .
\]  

However, $\bar{\mathcal{A}}$ is also a semi-free DGA in its own right.  The generators of $\A$ are naturally grouped into $\zp$ orbits, and if $\{x_i\}_{i=1}^n$ is a set containing exactly one representative from each orbit, then
\[
\bar{\mathcal{A}}(\widetilde{\Gamma})=T(\bar{x}_1, \bar{x}_2, ...\bar{x}_n\}.
\]

\begin{lemma}\label{lem:chain} $\phi$ is a chain map:
\[ \widetilde{\partial} \circ \phi =\phi \circ \widetilde{\partial}.
\]
\end{lemma}

Using Lemma~\ref{lem:chain}, we may define $\bar{\partial}: \bar{\mathcal{A}}(\widetilde{\Gamma})\rightarrow \bar{\mathcal{A}}(\widetilde{\Gamma})$ by
 \[
\bar{\partial}\bar{w}=\bar{\partial}(\phi(w))= \phi(\widetilde{\partial} w).
\]
This identifies the image of $\phi$ with the semi-free DGA   $(\bar{\mathcal{A}}(\widetilde{\Gamma}), \bar{\partial})$.

\begin{lemma}\label{lem:mainpf}
 $(\mathcal{A}(\Gamma), \partial)$ and $(\bar{\mathcal{A}}(\widetilde{\Gamma}), \bar{\partial})$ are isomorphic as semi-free DGAs.
\end{lemma}

\begin{lemma}\label{lem:equiv}
A $\zp$ stabilization of $(\A(\widetilde{\Gamma}), \gamma, \widetilde{\partial})$ induces an ordinary stabilization of $(\bar{\mathcal{A}}(\widetilde{\Gamma}), \bar{\partial})$, and a $\zp$ tame automorphism of $(\A(\widetilde{\Gamma}), \gamma, \widetilde{\partial})$ induces an ordinary tame automorphism of $(\bar{\mathcal{A}}(\widetilde{\Gamma}),\bar{\partial})$.
\end{lemma}

Theorem~\ref{thm:main} asserts that the equivalence type of $(\mathcal{A}(\Gamma), \partial)$ is an invariant of $K$.  Lemma~\ref{lem:mainpf} replaces this with a statement about $(\bar{\mathcal{A}}(\widetilde{\Gamma}), \bar{\partial})$.  Assuming Theorem~\ref{thm:triple} holds, Lemma~\ref{lem:equiv} then completes the proof.

\begin{proof}[Proof of Lemma~\ref{lem:chain}]
Let $x$ and $y$ be two generators of $\A(\widetilde{\Gamma})$.

\begin{align}
(\partial \circ \phi) (xy)&=\partial(\bar{x}\bar{y})\notag \\
&= \partial[ (\sum_{i=1}^p \gamma^i x)(\sum_{j=1}^p \gamma^j y)]\notag \\
&=\sum_{i=1}^p \sum_{j=1}^p \partial [(\gamma^i x)(\gamma^j y)]\notag \\
&= \sum_{i=1}^p \sum_{j=1}^p (\partial \gamma^i x)(\gamma^j y)+ (\gamma^ix)(\partial \gamma^j y)\notag\\
&= \sum_{i=1}^p \sum_{j=1}^p (\gamma^i \partial x)(\gamma^j y)+ (\gamma^ix)(\gamma^j \partial y)\notag
\end{align}

On the other hand, 
\begin{align}
(\phi \circ \partial )(xy)&= \phi  [(\partial x)y + x(\partial y)]\notag\\
&=\phi(\partial x) \phi (y) + \phi(x)\phi(\partial y))\notag\\
&= (\sum_{i=1}^p \gamma^i \partial x)(\sum_{j=1}^p \gamma^j y)+ (\sum_{i=1}^p \gamma^ix)(\sum_{j=1}^p \gamma^j \partial y)\notag
\end{align}
\end{proof}

\begin{proof}[Proof of Lemma~\ref{lem:mainpf}]

Each crossing in $\Gamma$ lifts to a $\zp$ orbit of crossings in $\widetilde{\Gamma}$.  If $x$ is a generator of $\A(\widetilde{\Gamma})$, there is a one-to-one correspondence between generators $\bar{x}\in \bar{\mathcal{A}}(\widetilde{\Gamma})$ and generators $\pi (x) \in \mathcal{A}(\Gamma)$.   The gradings of these generators agree, so $\bar{\mathcal{A}}(\widetilde{\Gamma})$  and $\mathcal{A}(\Gamma)$ are isomorphic as graded algebras.  

The boundary map on $\mathcal{A}(\Gamma)$ counts discs whose boundary has winding number $p$ with respect to the poles, and any such curve lifts to a  simple closed curve in $\widetilde{\Gamma}$.  On the other hand, any disc in $(S^2, \widetilde{\Gamma})$ projects to a disc in $(S^2, \Gamma)$ whose boundary has winding number $p$ with respect to the poles of $S^2$.  

By construction, the $x$ defect  of any admissible disc mapping into $(S^2, \widetilde{\Gamma})$ will agree with the $\pi(x)$ defect of its image in $\pi_*(S^2, \widetilde{\Gamma})$.  Thus a word $w(f,x)$ appears in $\partial \pi(x)$ if and only if $\phi(w(f, x))$ appears in the boundary of $\bar{x}$ in $\bar{\mathcal{A}}(\widetilde{\Gamma})$.

\end{proof}

\begin{proof}[Proof of Lemma~\ref{lem:equiv}]
The proof is almost immediate from the definitions.  A $\zp$ stabilization adds the generators $\bar{e}_1$ and $\bar{e}_2$ to $\bar{\mathcal{A}}(\widetilde{\Gamma})$, where $\bar{\partial} \bar{e}_1=\bar{e}_2$.  Similarly, if $f$ is an elementary $\zp$ automorphism which sends $a_i$ to $a_i+w_i$, then the map $\bar{f}$ which sends $\bar{f}(\bar{a})$ to $\bar{a}+\bar{w}$ is a tame automorphism of $\bar{\mathcal{A}}(\widetilde{\Gamma})$ which intertwines $\phi$.
\end{proof}

\section{Proof of Theorem~\ref{thm:triple}}\label{sect:triplepf}
Theorem~\ref{thm:triple} asserts that the equivalence type of the equivariant DGA $(\A(\widetilde{\Gamma}), \gamma, \widetilde{\partial})$ is an invariant of the Legendrian type of the knot $K \subset L(p,q)$.  This requires proving that Legendrian isotopy of $K$ changes $(\A(\widetilde{\Gamma}), \gamma, \widetilde{\partial})$ only by free $\zp$ stabilizations and $\zp$ tame isomorphisms.  When an isotopy occurs in the complement of the core curves $r_1=0$ and $r_2=0$, the proof is similar  to the proof of ordinary invariance for $(\A(\widetilde{\Gamma}), \widetilde{\partial})$.  However, the core curves are Reeb orbits where the bundle fails to be smooth, and more care must be taken with isotopies which pass $K$ across these fibers.  

\subsection{Reidemeister moves and isotopy away from the cores}

\begin{figure}[h]
\begin{center}
\scalebox{.35}{\includegraphics{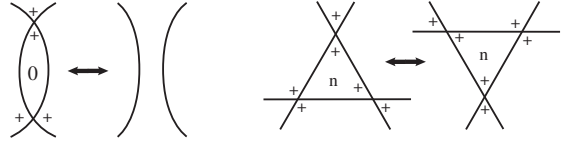}}
\end{center}
\caption{The two Reidemeister moves for Lagrangian projections.  In the right figure, $n\in \{0,1\}$.}\label{fig:Reids}
\end{figure}

\begin{lemma}[\cite{S}, Lemma 6.3]\label{lem:sablem}
If $\widetilde{\Gamma}_1$ and $\widetilde{\Gamma}_2$ are Lagrangian projections for Legendrian isotopic knots in $S^3$, then they differ by a sequence of the Reidemeister moves shown in Figure~\ref{fig:Reids}.
\end{lemma}

\begin{lemma} If $K_1$ and $K_2$ are Legendrian knots in $L(p,q)$ which are Legendrian isotopic in the complement of the core curves, then the Lagrangian projections $\widetilde{\Gamma}_1$ and $\widetilde{\Gamma}_2$ of their freely-periodic lifts differ by a sequence of $p$-tuples of the Reidemeister moves shown in Figure~\ref{fig:Reids}.
\end{lemma}

\begin{proof}   Away from the poles, the pair $(S^2, \widetilde{\Gamma})$ is a $p$-fold cover of the pair $(S^2, \Gamma)$.  Employing the argument in the proof of Lemma 6.3 of \cite{S}, the Lagrangian image of the isotopy is a homotopy of immersions away from the Reidemeister moves shown, and each of these  lifts to $p$ disjoint copies of the same move in $(S^2, \widetilde{\Gamma})$.
\end{proof}

\begin{proposition} If  $K_1$ and $K_2$ are Legendrian knots in $L(p,q)$ whose Lagrangian projections differ  by a Reidemeister move, then the equivariant DGAs $(\widetilde{\mathcal{A}}(\widetilde{\Gamma}_1), \gamma, \widetilde{\partial}_1)$ and $(\widetilde{\mathcal{A}}(\widetilde{\Gamma}_2), \gamma, \widetilde{\partial}_2)$  are equivalent. 
\end{proposition}

This proposition is proved in Section~\ref{sect:Reids}.

\subsection{Star moves and isotopy across the cores}\label{sect:stars}
\begin{figure}[h]
\begin{center}
\scalebox{.3}{\includegraphics{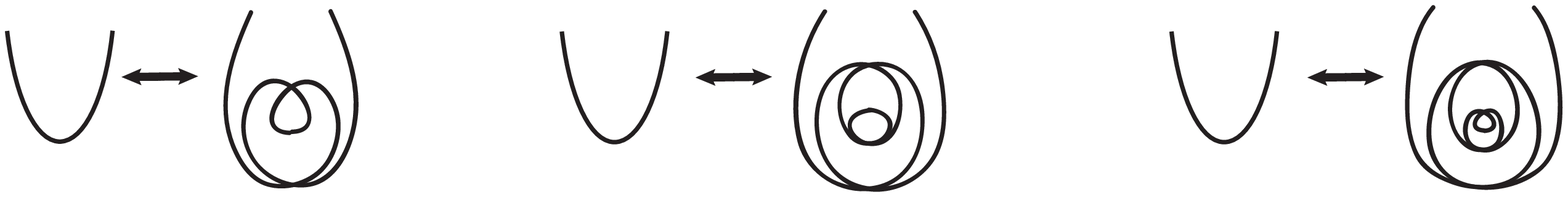}}
\end{center}
\caption{The change to $\Gamma$ induced by passing $K$ across a core curve. ($p=3,4,$ and $5$).}\label{fig:lpqstar}
\end{figure}

We turn now to isotopies which pass $K$ through the core of one of the Heegaard tori.  In a smooth bundle over the two-sphere, (e.g. $S^3$ or $L(p,1)$), such an isotopy would project to an isotopy of $\Gamma$ across one of the poles.  Recall, however, that $L(p,q)$ is the quotient space of $S^3$ under a cyclic group action which maps each core curve onto itself and each non-core curve into an orbit of $p$ fibers.  This implies the existence of  pairs of fibers in $L(p,q)$ which are not homotopic in the Heegaard tori: any fiber  on the boundary of an $\epsilon$ neighborhood of the core fiber has slope $\frac{x}{p}$  (where $x$ depends on the basis).  Thus there is in an orbifold point at each pole of the Lagrangian projection.  An isotopy passing $K$ across a core curve therefore passes the Lagrangian projection through a non-Reidemeister singularity, as shown in Figure~\ref{fig:lpqstar}. 

 \begin{figure}[h]
\begin{center}
\scalebox{.4}{\includegraphics{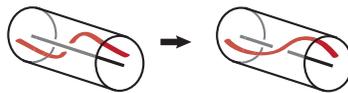}}
\end{center}
\caption{An isotopy passing $K$ across a core curve.}\label{fig:core}
\end{figure}

\begin{figure}[h]
\begin{center}
\scalebox{.3}{\includegraphics{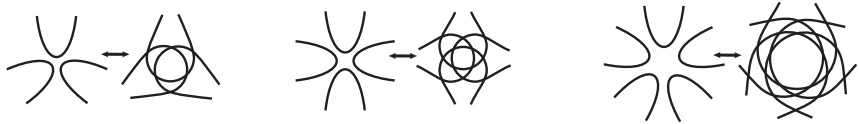}}
\end{center}
\caption{Star moves for $p=3,4,5$.  For each $p$, label the left-hand figure as $\widetilde{\Gamma}^-$ and the right-hand figure as $\widetilde{\Gamma}^+$.}\label{fig:sp}
\end{figure}

An isotopy passing $K$ across a core curve lifts to an isotopy passing $\K$ across a core curve in $S^3$ $p$ times.  The projection of this isotopy simultaneously moves $p$ strands of $\widetilde{\Gamma}$ across the corresponding pole of $S^2$ as shown in Figure~\ref{fig:sp}.

An isotopy passing $K$ across a core curve preserves the Legendrian type of the lift $\K$, so Lemma~\ref{lem:sablem} implies that $\widetilde{\Gamma}^-$ and $\widetilde{\Gamma}^+$  are related by a sequence of Reidemeister moves and the associated DGAs are stably tame isomorphic.  However, it is not clear that such a sequence respects the $\zp$ action.  This prompts the introduction of a new move relating generic Lagrangian diagrams.

\begin{definition} If $K_-$ and $K_+$ are Legendrian knots in $L(p,q)$ which differ only by an isotopy passing one strand across a core curve, then we say that a \textit{star move} relates the Lagrangian projections $\widetilde{\Gamma}^-$ and $\widetilde{\Gamma}^+$ of their freely periodic lifts.
\end{definition}

\begin{proposition}\label{prop:star}
 If two labeled diagrams $\widetilde{\Gamma}^-$ and $\widetilde{\Gamma}^+$ differ only by a star move, the equivariant DGAs $(\A(\widetilde{\Gamma}^-), \gamma, \widetilde{\partial}^-)$ and $(\A(\widetilde{\Gamma}^+), \gamma, \widetilde{\partial}^+)$ are $\zp$ equivalent.
\end{proposition}

The proof of this proposition will occupy Section~\ref{sect:star}

\noindent\textbf{Remark} Readers familiar with grid diagrams may note that there is a set of Legendrian grid moves which relates the grid diagrams of Legendrian equivalent knots in both $L(p,q)$ and $S^3$ \cite{BG}.  A single Legendrian grid move on a toroidal diagram in $L(p,q)$ corresponds to $p$ copies of a Legendrian grid move on a toroidal diagram in $S^3$.  This proves Proposition~\ref{prop:easyinvt}, but it is not strong enough to show Theorem~\ref{thm:triple}, as a grid move does not clearly translate into a $\zp$ equivalence of equivariant DGAs.

\subsection{Invariance under star moves}\label{sect:star}
This section is devoted to proving that a star move preserves the $\zp$ equivalence type of the equivariant DGA $(\A, \gamma, \widetilde{\partial})$. The proof is given for the case when $p$ is odd, but the proof for even $p$ is similar.  The structure of this argument is based on Chekanov's proof of Reidemeister II invariance in \cite{C}.  

 \begin{figure}[h]
\begin{center}
\scalebox{.5}{\includegraphics{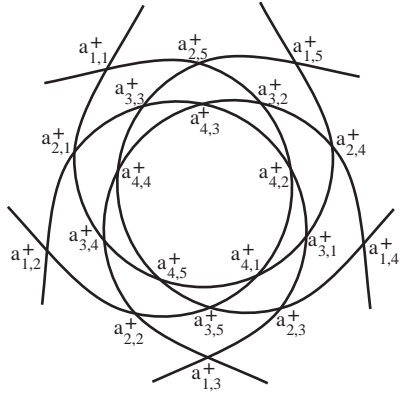}}
\end{center}
\caption{The star in $\widetilde{\Gamma}^+$, $p=5$.}\label{fig:splab}
\end{figure}

\subsubsection{A closer look at $(\mathcal{A}^+,\partial^+)$}

For simplicity, let $(\mathcal{A}^+,\partial^+)$ and $(\mathcal{A}^-, \partial^-)$ denote $(\A(\widetilde{\Gamma}^+), \gamma, \widetilde{\partial}^+)$ and $(\A(\widetilde{\Gamma}^-), \gamma, \widetilde{\partial}^-)$, respectively. Away from the star, assume all labels on the two diagrams agree. 

Consider first the algebra associated to $\widetilde{\Gamma}^+$, focusing only on the star region where it differs from $\widetilde{\Gamma}^-$.   The crossings in the star  are grouped into $p-1$ orbits so that  $\gamma(a_{j,i})=a_{j, i+1}$.  (See Figure~\ref{fig:splab}.) Furthermore, each crossing in $\widetilde{\Gamma}^+$ forms one corner of a bigon in the star, and the labels are assigned so that $a^+_{j,i}$ and $b^-_{p-j, i}$ appear at opposite ends of the same bigon.  

Recall the length function $\mathit{l}$ defined on generators of $\mathcal{A}^+$ in Section~\ref{sect:sab}, and note that it is constant on the generators in a fixed orbit.  In fact, the lengths of the orbits associated to the star crossings in $\widetilde{\Gamma}^+$ are clustered  around the values $\frac{m}{p}$, $m \in \mathbb{N}$.  The next lemma makes this statement precise and shows that the only maps involving multiple generators within such a cluster are given by the obvious bigons connecting $a_{k,i}$/$b_{p-k,i}$ pairs.  

 \begin{lemma}\label{lem:abndry}
If $b_{p-l,j}$ appears in a word $w \in \widetilde{\partial}a_{k,i}$ and $|\textit{l}(b_{p-l,j})-\mathit{l}(a_{k,i})| < \frac{1}{2p}$, then $w= b_{p-k, i}$.
\end{lemma}

\begin{proof}
The star results from passing $p$ segments of $\K$ across a core curve of $S^3$.  This isotopy can be assumed to take place in $p$ arbitrarily small balls distributed evenly around the core, and crossings in the star correspond to chords connecting the curve segments in different balls.  Thus the length associated to any generator is approximately $\frac{m}{p}$, $m\in \{ 1, 2, ...p-1 \}$, with the diameter of the ball providing an upper bound for the error.   

 \begin{figure}[h]
\begin{center}
\scalebox{.4}{\includegraphics{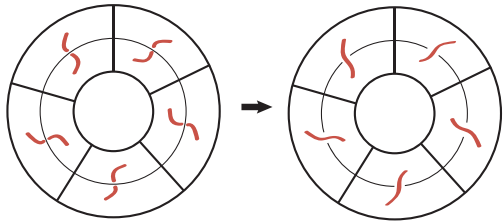}}
\end{center}
\caption{An isotopy of $\K$ which induces a star move creates generators whose lengths are roughly $\frac{m}{p}$.}\label{fig:chlength}
\end{figure}

Let $\epsilon'$ be the shortest length associated to a generator occurring outside the star, and set $\epsilon=\min(\epsilon', \frac{1}{4p})$.  Choosing the isotopy balls small enough ensures that the lengths of the star generators lie in the intervals $(\frac{m}{p}-\epsilon, \frac{m}{p}+\epsilon)$.  

Recall that an admissible disc $f:\Sigma\rightarrow S^2$ represents a term in the boundary of $a_{k,i}$ if and only if it satisfies
\[
0=\tilde{n}_{a_{k,i}}(f)=\frac{1}{2 \pi}\int_{\Sigma} f^*\Omega +\mathit{l}(a_{k,i}) - \sum_{y_j^- \in w(f,a_{k,i})} \mathit{l}(y_j).
\]  

We next show that $b_{p-k,i}$ always appears in $\partial a_{k,i}$.  Let $f$ be an admissible disc whose image is a bigon in the star.  The integral term in the defining equation for the defect can be made arbitrarily small by restricting the star to live in a sufficiently small neighborhood of the pole.  Since the length of each generator lies strictly between zero and one and the defect is always an integer, $\tilde{n}_{a_{k,i}}(f)$ must be zero.  Thus $b_{p-k,i}$ is as a summand in  $\widetilde{\partial} a_{k,i}$.

Finally, suppose that for some $m'$, both $\mathit{l}(a_{k,i})$ and $\mathit{l}(b_{p-l,j})$ lie in $(\frac{m'}{p}-\epsilon, \frac{m'}{p}+\epsilon)$.  The value $\mathit{l}(a_{k,i}) -  \mathit{l}(b_{p-l, j})$ is smaller than the length of any generator, so any word $w(a_{k,i}, f)\in \widetilde{\partial}a_{k,i}$ containing $b_{p-l, j}$ cannot contain any other generators.  This implies that the image of $f$ is a bigon, so $b_{p-l, j}=b_{p-k, i}$.
 \end{proof}

\subsubsection{Overview of proof}

In order to prove that  $(\mathcal{A}^+, \partial^+)$ and $(\mathcal{A}^-, \partial^-)$ are $\zp$ equivalent, we replace $\mathcal{A}^-$ with its free $\zp$ stabilization $\mathcal{A}'$:
\[ \mathcal{A}'=\mathcal{A}^- \coprod \mathcal{E}^p_1\coprod \mathcal{E}^p_2 ...\coprod \mathcal{E}^p_{p-1}. \] 
The $2(p-1)$ new orbits are assigned gradings so that  $(\mathcal{A}', \partial')$ is isomorphic to $(\mathcal{A}^+, \partial^+)$ as a graded algebra:
\[|e^k_{1,i}|=|a_{k,i}|\text{ for  }1\leq k \leq p-1,\]
 and in fact they are tamely isomorphic equivariant DGAs.  The proof begins with the construction of an explicit isomorphism $s:\mathcal{A}^+\rightarrow \mathcal{A}'$ (Section~\ref{sect:s}).  Conjugating by $s$ gives a new boundary map $\hat{\partial}$ on $\mathcal{A}'$:
\[
\hat{\partial}=s\circ \partial^+ \circ s^{-1}.
\]

 \[
\begin{diagram}
\node{\mathcal{A}^+}\arrow{e,t}{s}
\node{ \mathcal{A}'} 
\arrow{s,r}{g}
\node{\mathcal{A}^-}
\arrow{w,t}{\zp \text{ stab.}}\\
\node{}
\arrow{e,t,!}{}
\node{\mathcal{A}'}
\end{diagram}
\]
 
 In Section~\ref{sect:g} we  define a graded $\zp$ tame isomorphism  $g:\mathcal{A}'\rightarrow\mathcal{A}'$.  The map $g$ is constructed so that $\partial'=g\circ\hat{\partial} \circ g^{-1}$, which implies the $\zp$ tame isomorphism of $\mathcal{A}^+$ and $\mathcal{A}'$.

\subsubsection{The map $s:\mathcal{A}^+\rightarrow \mathcal{A}'$}\label{sect:s}
As a first step towards defining the isomorphism $s:\mathcal{A}^+\rightarrow \mathcal{A}'$, we extend the length function $\mathit{l}$ to generators of $\mathcal{A}'$.  The generators of $\mathcal{A}'$ are of two types: generators coming from crossings in $\widetilde{\Gamma}^-$ and generators in the $\mathcal{E}_k^p$. To a generator of the first type, assign the length of the corresponding generator in $\mathcal{A}^-$.  Similarly, the star generators of $\mathcal{A}^+$ can be used to assign lengths to the generators of $\mathcal{A}'$ coming from the stabilizing orbits.  Set
\begin{align}
\mathit{l}(e_{1,i}^k)&=\mathit{l}(a_{k,i})\notag \\
\mathit{l}(e_{2,i}^k)&=\mathit{l}(b_{p-k,i}).\notag
\end{align}

Note that this implies $\mathit{l}(e_{j,i}^k)$ lies in $(\frac{m}{p}-\epsilon, \frac{m}{p}+\epsilon)$ for some $m\in \{1,2, ...p-1\}$.  Perturbing $K$ slightly and letting $\epsilon\rightarrow 0$, we may assume that no other generators' lengths lie in these intervals.  

Let $\mathcal{O}_{[0]}$ denote the set of generators of $\mathcal{A}'$ whose length is less than or equal to $\frac{1}{p}+\epsilon$.  
Order the remaining orbits by increasing length, and label them as $\mathcal{O}_{[j]}, j=1, 2, 3...$  so that if $x_j \in \mathcal{O}_{[j]},$ then $\mathit{l}'(x_{1})<\mathit{l}'(x_2)<\mathit{l}'(x_3)...$.  Let $\mathcal{A}_{[j]}$ be the subalgebra generated by the elements $ \mathcal{O}_{[i]}$ for $0\leq i \leq j$.
\[ 
A_{[j]} = T(\mathcal{O}_{[0]}, \mathcal{O}_{[1]}, ....\mathcal{O}_{[j]}).
\]

\begin{lemma}\label{lem:length} If $x_k$ is in $\mathcal{O}_{[j]}$ for $j\geq1$, then $\partial'(x_k) \subset \mathcal{A}_{[j-1]}$.
\end{lemma}

\begin{proof} If $w(f,x_l)$ appears in $\partial'x_k$, then $\tilde{n}_{x_k}(f)=0$: 
\[
0=\frac{1}{2\pi}\int_{\Sigma} f^*\Omega +\mathit{l}(x_k) - \sum_{x_l^- \in w(f,x_k)} \mathit{l}(x_l).
\]
The integral term is negative, and each $\mathit{l}(x_l)$ is positive, which proves the lemma.
\end{proof}

For each $a_{k,i}$, write $\partial^+a_{k,i}=b_{p-k,i}+v_{k,i}+w_{k,i}$, where words in $v_{k,i}$ involve only generators from crossings outside the star, and each word in $w_{k,i}$ contains at least one generators coming from a star crossing.   Let $\mathcal{O}_M = \{ x\in \mathcal{A}^+ | \mathit{l}'(x)\leq \frac{M}{p}+\epsilon\}$.  If $a_{k,i} \in \mathcal{O}_M \backslash \mathcal{O}_{M-1}$, then  Lemma~\ref{lem:abndry} shows that $w_{k,i} \in \mathcal{O}_{M-1}$.  We define  maps $s_M:\mathcal{O}_M \rightarrow \mathcal{A}'$ inductively: 
\[
s_1(x) = \begin{cases}  
 x& \text{if }x \text{ comes from a crossing in } \widetilde{\Gamma}^-\\
e_{1,i}^k &  \text{if} \ x=a_{k,i} \\ 
e_{2,i}^k + v_{k,i} & \text{if}\ x=b_{p-k,i}. \end{cases} 
\]
For the inductive step, suppose that $s_i$ is defined for $i\in \{ 1,2, ...M-1\}$.
\[
s_M(x) = \begin{cases}  
 x& \text{if }x \text{ comes from a crossing in } \widetilde{\Gamma}^-\\
e_{1,i}^k &  \text{if} \ x=a_{k,i} \\ 
e_{2,i}^k + v_{k,i} + s_{M-1}(w_{k,i}) & \text{if}\ x=b_{p-k,i}. \end{cases} 
\]
\begin{definition} 
Define $s:\mathcal{A}^+ \rightarrow \mathcal{A}'$  by 
\[
s(x) = \begin{cases}  
 x&\text{if }x \text{ comes from a crossing in } \widetilde{\Gamma}^-\\
e_{1,i}^k &  \text{if} \ x=a_{k,i} \\ 
e_{2,i}^k + v_{k,i} +s_{p-1}(w_{k,i})& \text{if}\ x=b_{p-k,i}. \end{cases} 
\]
\end{definition}

By construction, $s$ preserves gradings and intertwines the $\zp$ actions on $\mathcal{A}^+$ and $\mathcal{A}'$.

\subsubsection{The projection map $\tau:\mathcal{A}'\rightarrow \mathcal{A}^-$}\label{sect:tau}
As a vector space, $\mathcal{A}'$ decomposes as $\mathcal{A}^-\oplus \mathcal{I}_{\mathcal{E}}$, where $\mathcal{I}_{\mathcal{E}}$ is the two-sided ideal generated by elements in the $\mathcal{E}^p_k$.  Define $F: \mathcal{A}'\rightarrow \mathcal{A}'$ by
\[
F(x) = \begin{cases}  ye_{1,i}^kz &  \text{if} \ x=ye_{2,i}^kz \text{ and }y\in \mathcal{A}^-\\ 
0 & otherwise.\end{cases} 
\]
Note that $F$ is graded with degree $1$.

Let $\tau: \mathcal{A}'\rightarrow \mathcal{A}'$ be projection to $\mathcal{A}^-$.
\begin{lemma}
\ $\tau$ satisfies
\begin{equation}\label{lem:tau}
\tau +id_{\mathcal{A}'}=F \circ \partial ' + \partial' \circ F.
\end{equation}
\end{lemma} 
The proof is a straightforward computation.

\begin{lemma}\label{lem:taucomp}
$\tau \circ \partial'=\tau \circ \hat{\partial}.$
\end{lemma}

\begin{proof}
We first show that on the $\mathcal{E}_k^p$, the stronger statement $\partial'=\hat{\partial}$ holds.  This fact will be used of this in Section~\ref{sect:g}.

First recall that $\partial'(e_{1,i}^k)=e_{2,i}^k$ and $\partial'(e_{2,i}^k)=0$.  Compare this to $\hat{\partial}(e_{j,k}^k)$:
\begin{align}
\hat{\partial}(e_{1,i}^k)&=s\circ \partial^+ \circ s^{-1}(e_{1,i}^k)\notag\\
&=s\circ \partial^+(a_{k,i})\notag\\
&= s (b_{p-k,i}+v_{k,i}+w_{k,i})\notag\\
&= (e_{2,i}^k+v_{k,i}+s_{p-1}(w_{k,i}))+s(v_{k,i})+s(w_{k,i})\notag\\
&=e_{2,i}^k\notag\\
\newline\notag\\
\hat{\partial}(e_{2,i}^k)&=s\circ \partial^+ \circ s^{-1}(e_{2,i}^k)\notag\\
&=s\circ \partial^+(b_{p-k,i}+ v_{k,i}+w_{k,i})\notag\\
&=s\circ \partial^+(\partial^+a_{k,i})\notag\\
&=0\notag
\end{align}

Now consider some generator $x$ which is associated to a crossing in $\widetilde{\Gamma}^-$.  To prove the lemma, we compare the words appearing in $\tau \circ \partial'(x)$ and in $\tau \circ \hat{\partial}(x)$; it suffices to show that any word with no generators in the $\mathcal{E}^p_k$ appears in both $\hat{\partial}(x)$ and $\partial'(x)$.  This argument is similar to the proof of Step 5 in \cite{S}.

Terms in $\partial'x$ come from one of the following types of discs  (Figure~\ref{fig:neck}):
\begin{enumerate}
\item Discs which have the same multiplicity throughout the star region;
\item Discs which flow through the star region in $\widetilde{\Gamma}^-$.
\end{enumerate} 
Clearly, $\tau\circ \partial'(x)=\partial'(x)$.

 \begin{figure}[h]
\begin{center}
\scalebox{.4}{\includegraphics{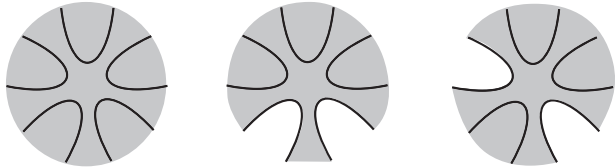}}
\end{center}
\caption{Discs in $\partial'$.  Left:  A disc with constant multiplicity in the star. Center and Right: Two discs flowing through the star.}\label{fig:neck}
\end{figure}

Now expand $\hat{\partial}$ as $s\circ \partial^+ \circ s^{-1}$.
Since $s^{-1}(x)=x$, terms in $\hat{\partial}x$ are the $s$-images of terms in $\partial^+$.  These come in three flavors:
\begin{enumerate}
\item\label{1a} Words involving none of the generators associated to crossings in the star;
\item\label{2} Words involving some $a_{k,i}$ but no $b_{j,l}$;
\item\label{3} Words involving $b_{j,i}$.
\end{enumerate} 

Boundary terms of the first type in both lists agree.  In the second list, note that $s$ sends words involving only $a_{k,i}$ terms to $\mathcal{I}_{\mathcal{E}}$, so these will vanish under $\tau$.  

To see that the remaining terms agree, recall that  $s(b_{p-k,i})=e_{2,i}^k+v_{k,i}+s_{p-1}(w_{k,i})$.  Any word containing $e_{2,i}^k$ vanishes under  $\tau$. However, the terms coming from $v_{k,i}+s_{p-1}(w_{k,i})$ which  are not killed by $\tau$ represent boundary discs which start at  $a_{k,i}^+$.  These can be ``glued" to boundary discs for $x$ with a corner at $b_{p-k,i}^-$ to produce boundary discs for $x$ in $\widetilde{\Gamma}^-$. See Figure~\ref{fig:glue}. However, these are exactly the discs in $\partial'x$ which flow through the star.  The proof that this gluing operation is smooth comes from the argument in  \cite{S}.

 \begin{figure}[h]
\begin{center}
\scalebox{.4}{\includegraphics{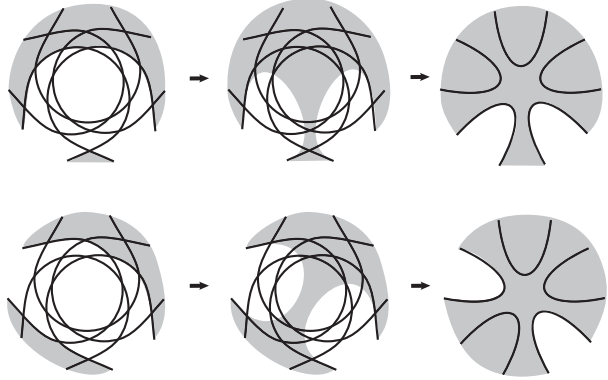}}
\end{center}
\caption{Gluing discs in $\partial^+x$ with which hit $b_{p-k,i}^-$ to discs in $\partial^+a_{k,i}^+$ gives discs in $\partial'x$.  Top: $b_{1,3}$. Bottom: $b_{2,2}$.}\label{fig:glue}
\end{figure}

\end{proof}

\subsubsection{Constructing $g:\mathcal{A}'\rightarrow \mathcal{A}'$}\label{sect:g}

Following Chekanov, we construct $g:\mathcal{A}'\rightarrow \mathcal{A}'$ as a composition of maps $g_j$.  Each $g_j$ is a graded $\zp$ elementary isomorphism which is the identity away from the orbit  $\mathcal{O}_j$. Furthermore, each  $g_j$ inductively defines a new boundary map $\partial_{[j]}$ on $\mathcal{A}'$  by conjugation:
\[
\partial_{[j]}=g_{j}\partial_{[j-1]}g_{j}^{-1}.
\]

Setting $g=g_n \circ g_{n-1}\circ ...\circ g_2 \circ g_1$, we will prove that 
\[
\partial'=g\circ \hat{\partial} \circ g^{-1}=g\circ s\circ \partial^+\circ s^{-1}\circ g^{-1}.
\]
This establishes a $\zp$ tame isomorphism between $(\mathcal{A}^+,\partial^+)$ and $(\mathcal{A}', \partial')$, and thus the $\zp$ equivalence of $(\mathcal{A}^+, \partial^+)$ and $(\mathcal{A}^-, \partial^-)$.

The $g_j$ should satisfy $\partial_{[j]}|_{\mathcal{A}_{[j]}}=\partial'|_{\mathcal{A}_{[j]}}$.  We begin by setting $\partial_{[0]}=\hat{\partial}$.  As noted in the proof of Lemma~\ref{lem:taucomp},  $\hat{\partial}=\partial'$ on the generators in the stabilizing orbits, and Lemma~\ref{lem:length} implies that they also agree on generators with length less that $\frac{1}{p}-\epsilon$.  This establishes the base case $\partial_{[0]}|_{\mathcal{A}_{[0]}}=\partial'|_{\mathcal{A}_{[0]}}$.   

For the inductive step, suppose that for $1\leq k \leq j-1$, the maps $g_k$ satisfy $\partial_{[k]}|_{A_{[k]}}=\partial'|_{A_{[k]}}$.  Define $g_j: \mathcal{A}'\rightarrow\mathcal{A}'$ by
\[
g_j(x) = \begin{cases} x+F(\partial'(x)+\partial_{[j-1]}(x)) & \text{if}\ x \in \mathcal{O}_j\\ 
x & \text{otherwise}. \end{cases} 
\]

It follows from Lemma~\ref{lem:length}, the proof of Lemma~\ref{lem:taucomp}, and the definition of $F$ that if $x\in \mathcal{O}_{[j]}$, then $F(\partial'(x)+\partial_{[j-1]}(x)) \in A_{[j-1]}$.   
By the inductive hypothesis, the restriction of $\partial_{[j-1]}$ to $A_{[j-1]}$ agrees with $\partial'$.  We have the following for $x \in \mathcal{O}_{[j]}$:
\begin{align}
\partial_{[j]}(x)& =g_j \circ \partial_{[j-1]} \circ g_j^{-1}(x) \notag \\
 &=g_j \circ \partial_{[j-1]}(x+F(\partial'(x)+\partial_{[j-1]}(x))) \notag \\
     &= g_j \circ \partial_{[j-1]}(x) +g_j\circ \partial_{[j-1]}(F(\partial'(x)+\partial_{[j-1]}(x)))\notag \\
    &= \partial_{[j-1]}(x) +g_j \circ \partial'(F(\partial'(x)+\partial_{[j-1]}(x))) \notag \\
  &=  \partial_{[j-1]}(x) +(\partial'\circ F)(\partial'(x)+\partial_{[j-1]}(x))\notag 
\end{align}

Again taking $x\in \mathcal{O}_{[j]}$, we apply Lemma~\ref{lem:tau} to the last line:
 \begin{align}
\partial_{[j]}(x) &=  \partial_{[j-1]}(x)+(\tau+Id_{\mathcal{A}'}+F \circ \partial') (\partial'(x)+\partial_{[j-1]}(x))\notag \\
&= \partial_{[j-1]}(x)+ \tau\circ \partial'(x)+\tau\circ \partial_{[j-1]}(x)  + \partial'(x)+\partial_{[j-1]}(x)+  F\partial'\partial'(x)+F\partial'\partial_{[j-1]}(x).\notag \\
&= \tau\circ \partial'(x)+\tau\circ \partial_{[j-1]}(x) ) + \partial'(x) \notag 
 \end{align}

Expand the middle term in the last line:
 \[\tau\circ \partial_{[j-1]}(x)=\tau\circ g_{j-1}\circ g_{j-2}\circ...\circ g_1\circ \hat{\partial}(x)= \tau\circ \hat{\partial}(x).\]
 
 Thus \[\tau\partial'(x)+\tau\partial_{[j-1]}(x) =(\tau\circ \partial'+\tau\circ \hat{\partial})(x)=0\]
 by Lemma~\ref{lem:taucomp} and
 \[\partial_{[j]}(x)= \partial'(x)\]
 as desired. This proves that $\partial'$ and $\partial_{[j]}$ agree on $\mathcal{A}_{[j]}$.  Inductively, this implies $g\circ \hat{\partial} \circ g^{-1} =\partial'$, which completes the proof of Proposition~\ref{prop:star}.
 
\subsection{Reidemeister invariance}\label{sect:Reids}
To complete the proof of Theorem~\ref{thm:triple}, we show that the equivalence type of the equivariant DGA $(\A(\widetilde{\Gamma}), \gamma, \widetilde{\partial})$ is preserved by a Reidemeister move on $\Gamma$.   Recall that each Reidemeister move in $(S^2,\Gamma)$ lifts to a $p$-tuple of Reidemeister moves in $(S^2, \widetilde{\Gamma})$.  

\begin{figure}[h]
\begin{center}
\scalebox{.4}{\includegraphics{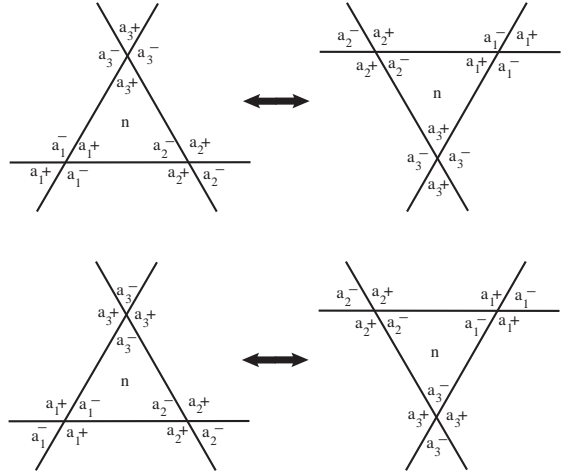}}
\end{center}
\caption{The two possible labelings for the Reidemeister III move, where $n\in \{0,1\}$.}\label{fig:r3}
\end{figure}

A Reidemeister II move on $\Gamma$ adds $2p$ new crossings to $\widetilde{\Gamma}$.  The proof that the ``before" and ``after" diagrams yield equivalent equivariant DGAs is similar to the proof of Proposition~\ref{prop:star}, so we turn to  Reidemeister III.  This argument is based on the proof in \cite{S}.  Up to rotation and switching $a$- and $b$-type generators, there are two versions of this move in $(S^2, \Gamma)$.   (See Figure~\ref{fig:r3}).  Each of these lifts to a $p$-tuple of identical local moves in $(S^2, \widetilde{\Gamma})$.  

Suppose that each of the Reidemeister triangles in $\widetilde{\Gamma}_1$ (respectively,  $\widetilde{\Gamma}_2$)  look locally like the left (right) diagram shown in the top row of Figure~\ref{fig:r3}. (This is the case left to the reader in \cite{S}.)  Set $n=0$. Since $\A(\widetilde{\Gamma}_1)$ and $\A(\widetilde{\Gamma}_2)$ have the same number of generators with the same labels, we identify the algebras and show that the boundary maps corresponding to the two diagrams give isomorphic equivariant DGAs.  

Define the  the tame isomorphisms $f_i: (\A, \widetilde{\partial}_1) \rightarrow (\A, \widetilde{\partial}_2)$:
\[
f_i(x) = \begin{cases}  a_{2,i}+a_{1,i}a_{3,i}&  \text{if} \ x=a_{2,i} \\ 
b_{1,i}+a_{3,i}b_{2,i} &\text{if} \ x=b_{1,i}\\
b_{3,i}+b_{2,i}a_{1,i} &\text{if}\ x=b_{3,i} \\
x &\text{otherwise.}  \end{cases} 
\]

Let $f=f_p\circ f_{p-1}\circ ...\circ f_1$.

\begin{lemma} The map $f$ is a graded tame automorphism which satisfies $f \circ \widetilde{\partial}_1=\widetilde{\partial}_2 \circ f$.
\end{lemma}

\begin{proof}  That $f$ is tame follows from the definition.  To see that  $f$ is graded, consider a generator $x$ associated to a crossing away from the Reidemeister triangle, and suppose that some disc representing a term in  $\widetilde{partial}_2x$ which crosses the triangle and has a corner at $a_2^-$.  Since $n=0$, truncating the disc at the $a_1a_3$ edge gives a new disc which has the same defect.  This amounts to replacing $a_2^-$ by $a_1^-a_3^-$  and getting a new boundary word.  Since the boundary map is graded, this implies that $|a_1a_3|=|a_2|$.  The arguments for the other generators are similar.

To prove $f \circ \widetilde{\partial}_1=\widetilde{\partial}_2 \circ f$, apply the two compositions to an  arbitrary generator and compare the resulting terms.  We demonstrate this comparison for $x=b_{1,i}$.

Figure~\ref{fig:bndry1} shows that $\partial(b_{1,i})=Ub_{2,i} + a_{3,i}V + W + X$, where each capital letter represents a sum of words not involving any of the other local generators.  The map $f_i$ fixes $b_{2,i}$ and $a_{3,i}$, so
\[
(f _i \circ \widetilde{\partial}_1 )(b_{1,i})=Ub_{2,i} + a_{3,i}V + W + X.
\]
\begin{figure}[h]
\begin{center}
\scalebox{.3}{\includegraphics{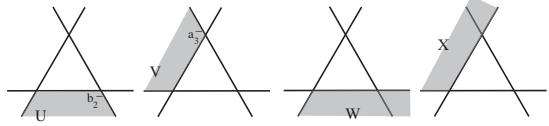}}
\end{center}
\caption{Terms in $\partial(b_{1,i})$.}\label{fig:bndry1}
\end{figure}

On the other hand, $f_i (b_{1,i})=b_{1,i}+a_{3,i}b_{2,i}$ and Figure~\ref{fig:bndry2} shows the following:
\begin{align}
( \widetilde{\partial}_2 \circ f _i )(b_{1,i})& = \widetilde{\partial}_2b_{1,i}+ (\widetilde{\partial}_2 a_{3,i})b_{2,i}+a_{3,i}(\partial b_{2,i})\notag \\
 &=W + a_{3,i}Z+ X +Yb_{2,i}+(U+Y)b_{2,i}+a_{3,i}(V+Z).\notag
\end{align}

Working modulo two, this shows $f \circ \widetilde{\partial}_1(b_{1,i})=\widetilde{\partial}_2 \circ f(b_{1,i})$.  The arguments are similar for generators associated to other crossings.  

\begin{figure}[h]
\begin{center}
\scalebox{.3}{\includegraphics{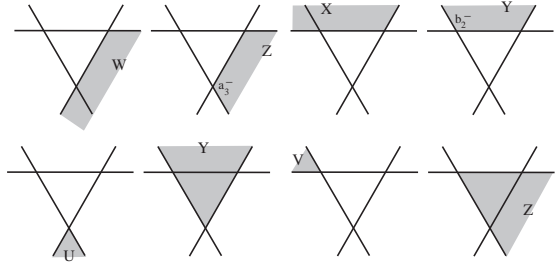}}
\end{center}
\caption{Terms in $\partial (b_{1,i}+a_{3,i}b_{2,i})$.}\label{fig:bndry2}
\end{figure}

It remains to show that the map $f$ commutes with the $\zp$ action $\gamma$. First note that $f\circ \gamma=f$ because the $p$ Reidemeister triangles are permuted by $\gamma$. 

Although $f$ is defined as $f_p\circ f_{p-1}\circ ...\circ f_1$,  in fact it is independent of the order of the composition.  This allows us to write $f=\sum_{i=1}^p f_i$.  Since $ \gamma \circ f_i =f_{i+1}$, we have 
\[
\gamma \circ f=\gamma\circ (f_1 +f_2 + ... + f_p)=f_ 2  ...+ f_p + f_1=f.
\] 

\end{proof}

This  completes the proof for the chosen case. When  $n=1$, the map which sends each generator of $\A(\widetilde{\Gamma}_1)$ to the generator of $\A(\widetilde{\Gamma}_2)$ with the same label intertwines the $\widetilde{\partial}_i$.  The proof for the other Reidemeister III move is given explicitly in \cite{S}, and the argument is similar to the one provided here.

\section{Examples}
We conclude with two examples.

\subsection{A pair of knots which are not Legendrian isotopic}
\begin{figure}[h]
\begin{center}
\scalebox{.5}{\includegraphics{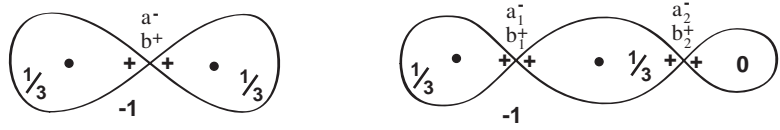}}
\end{center}
\caption{Left: $K_1$. RIght: $K_2$, which is a stabilization of $K_1$. The dots represent the poles of $S^2$. }\label{fig:ex2}
\end{figure}

For the first example, we consider two knots in the lens space $L(3,2)$.  (See Figure~\ref{fig:ex2}.)  The two knots differ by a Legendrian stabilization, a topological isotopy which does not preserve the Legendrian type of the knot.   

The boundary map $\partial_1$ is the zero map, so the homology of $(\mathcal{A}_1, \partial_1)$ is just a free group on two generators.  

In the case of the stabilized knot, the images of the generators are as follow:
\begin{align}
\partial_2 a_1&=1+b_2+a_2a_2\notag\\
\partial_2 b_1&= a_1a_2a_2 + a_2a_2a_1\notag\\
\partial_2 a_2&=1\notag\\
\partial_2 b_2&= a_2a_1a_1+a_1a_1a_2.
\end{align}

In general, distinguishing equivalence classes of DGAs can be difficult, and a variety of algebraic tools have been developed to make this problem more tractable.  We refer the reader to \cite{C} and \cite{Ng} for a discussion of  Chekanov polynomials, augmentations, linearized homology, and the characteristic algebra, but we note that the following suffices to distinguish $K_1$ and $K_2$:

\begin{proposition} An \textit{augmentation} of a DGA is an algebra homomorphism $\epsilon:\mathcal{A}\rightarrow \mathbb{Z}_2$ such that $\epsilon(1)=1$, $\epsilon \circ \partial=0$ and $\epsilon(x)=0$ if $|x|\neq 0$.  The existence of augmentations is an invariant of the equivalence type of the algebra.
\end{proposition}

The identity map is an augmentation of $(\mathcal{A}_1, \partial_1)$, whereas $(\mathcal{A}_2, \partial_2)$ has no augmentations.  Thus the two DGAs are not equivalent, and $K_1$ is not Legendrian isotopic to $K_2$. The knots in this example can also be distinguished by the classical invariants of their lifts to $S^3$; we would be interested in studying pairs of Legendrian non-isotopic knots in $L(p,q)$ which are not distinguished by classical invariants.

\subsection{An example in $L(5,2)$}
For the second example, we compute $(\mathcal{A}, \partial)$ for a knot in $L(5,2)$.  The Lagrangian projection is shown as a rectangle; to recover $S^2$, collapse each of the top and bottom edges to a point and identify the vertical edges as indicated.  The knot shown here is $K(5,2,2)$ in the notation of \cite{R2}.

\begin{figure}[h]
\begin{center}\notag
\scalebox{.6}{\includegraphics{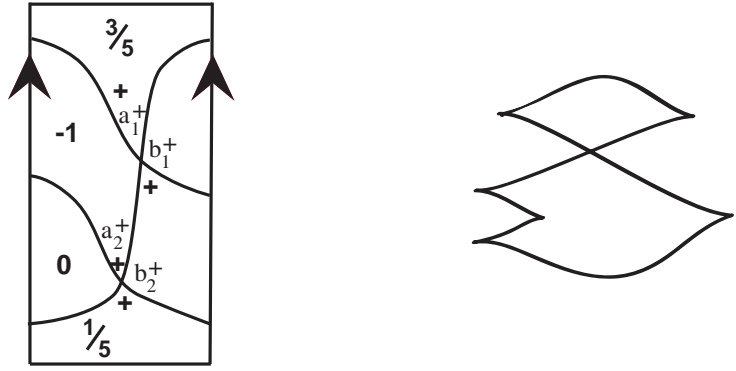}}
\end{center}
\caption{Left: A labeled diagram for $K\subset L(5,2)$. Right: A front projeciton for $\widetilde{K}\subset S^3$.}\label{fig:exdiag}
\end{figure}

The cyclic group grading $\mathcal{A}$ is $\mathbb{Z}_2$, and we have $|a_i|=1$ and  $|b_i|=0$ for $i=1,2$.
\begin{align}
\partial a_1&=a_2a_2\notag \\
\partial b_1&=0\notag \\
\partial a_2&=1\notag \\
\partial b_2&=a_2 b_1+b_1a_2\notag \\
\end{align}

\begin{figure}[h]
\begin{center}\notag
\scalebox{.6}{\includegraphics{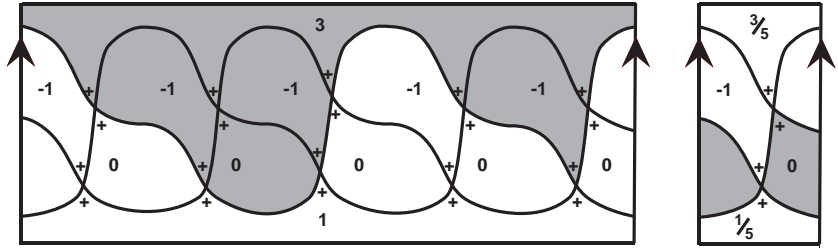}}
\end{center}
\caption{Discs which represent terms in $\partial$.  The disc bounded by a capping path for $a_2$ is shown lifted to $\widetilde{\Gamma}$.}\label{fig:exdown}
\end{figure}

In addition, both capping paths for $a_1$ bound admissible discs with $\tilde{n}_{a_2}(f)=0$; the corresponding terms cancel modulo two, and these discs are not shown. We note that although $\partial a_2=1$, $K$ is not a stabilization of any other Legendrian knot in $L(5,2)$.  The proof of this fact relies on the classification of Legendrian unknots in $S^3$ due to Eliashberg and Fraser \cite{EF}.

\begin{figure}[h]
\begin{center}
\scalebox{.3}{\includegraphics{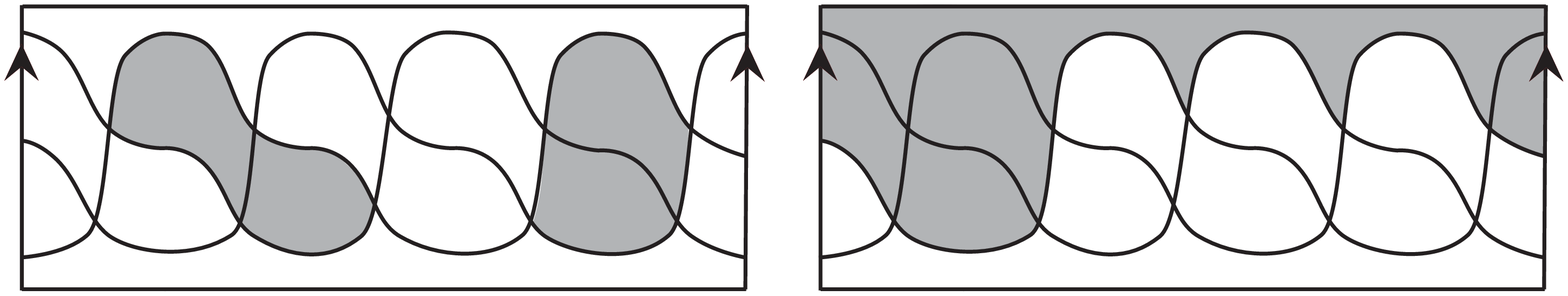}}
\end{center}
\caption{Discs representing terms in $\widetilde{\partial}$ which cancel modulo two under $\phi$.  The labels agree with the labels in the left-hand diagram in Figure~\ref{fig:exdown}.}\label{fig:ex}
\end{figure}

Lemma~\ref{lem:mainpf} implies that we could also compute the differential from the $\phi$-image of $\widetilde{\partial}$ in $\A(\widetilde{\Gamma})$.   In this computation there are additional discs representing non-canceling terms in $\widetilde{\partial}$ which nevertheless cancel in the image of $\phi$.  Representatives of these discs are shown in Figure~\ref{fig:ex}. In $(\A(\widetilde{\partial}), \widetilde{\partial})$, the terms in $\widetilde{\partial}(x_{i,0})$ are as follow:
\begin{align}
\widetilde{\partial} a_{1,0}&=b_{1,1}+a_{2,0}a_{2,1}+b_{1,-1}\notag\\
\widetilde{\partial} b_{1,0}&=a_{2,2}+a_{2,-1}\notag\\
\widetilde{\partial} a_{2,0}&=1\notag\\
\widetilde{\partial} b_{2,0}&=a_{1,-2}+a_{1,1}+a_{2,1}b_{1,0}+b_{1,-1}a_{2,-1}\notag
\end{align}

\bibliographystyle{alpha}
\bibliography{contactbib.bib}

\begin{thebibliography}{EGH00}

\bibitem[BG]{BG}
Kenneth~L. Baker and J.~Elisenda Grigsby.
\newblock Grid {D}iagrams and {L}egendrian {L}ens {S}pace {L}inks.
\newblock arXiv:0804.3048.

\bibitem[Che02]{C}
Yuri Chekanov.
\newblock Differential algebra of {L}egendrian links.
\newblock {\em Invent. Math.}, 150(3):441--483, 2002.

\bibitem[EF98]{EF}
Yakov Eliashberg and Maia Fraser.
\newblock Classification of topologically trivial {L}egendrian knots.
\newblock In {\em Geometry, topology, and dynamics ({M}ontreal, {PQ}, 1995)},
  volume~15 of {\em CRM Proc. Lecture Notes}, pages 17--51. Amer. Math. Soc.,
  Providence, RI, 1998.

\bibitem[EGH00]{EGH}
Y.~Eliashberg, A.~Givental, and H.~Hofer.
\newblock Introduction to symplectic field theory.
\newblock {\em Geom. Funct. Anal.}, (Special Volume, Part II):560--673, 2000.
\newblock GAFA 2000 (Tel Aviv, 1999).

\bibitem[Eli98]{E}
Yakov Eliashberg.
\newblock Invariants in contact topology.
\newblock In {\em Proceedings of the {I}nternational {C}ongress of
  {M}athematicians, {V}ol. {II} ({B}erlin, 1998)}, number Extra Vol. II, pages
  327--338 (electronic), 1998.

\bibitem[ENS02]{ENS}
John~B. Etnyre, Lenhard~L. Ng, and Joshua~M. Sabloff.
\newblock Invariants of {L}egendrian knots and coherent orientations.
\newblock {\em J. Symplectic Geom.}, 1(2):321--367, 2002.

\bibitem[Etn03]{Et}
John~B. Etnyre.
\newblock Introductory lectures on contact geometry.
\newblock In {\em Topology and geometry of manifolds ({A}thens, {GA}, 2001)},
  volume~71 of {\em Proc. Sympos. Pure Math.}, pages 81--107. Amer. Math. Soc.,
  Providence, RI, 2003.

\bibitem[Gei08]{Ge}
Hansj{\"o}rg Geiges.
\newblock {\em An introduction to contact topology}, volume 109 of {\em
  Cambridge Studies in Advanced Mathematics}.
\newblock Cambridge University Press, Cambridge, 2008.

\bibitem[GRS08]{GRS}
J.~Elisenda Grigsby, Daniel Ruberman, and Sa{\v{s}}o Strle.
\newblock Knot concordance and {H}eegaard {F}loer homology invariants in
  branched covers.
\newblock {\em Geom. Topol.}, 12(4):2249--2275, 2008.

\bibitem[HLN06]{HLN}
Jonathan~A. Hillman, Charles Livingston, and Swatee Naik.
\newblock Twisted {A}lexander polynomials of periodic knots.
\newblock {\em Algebr. Geom. Topol.}, 6:145--169 (electronic), 2006.

\bibitem[Ng03]{Ng}
Lenhard~L. Ng.
\newblock Computable {L}egendrian invariants.
\newblock {\em Topology}, 42(1):55--82, 2003.

\bibitem[Ras07]{R2}
J.~A. Rasmussen.
\newblock Lens space surgeries and {L}-space homology spheres.
\newblock arXiv:0710.2531, 2007.

\bibitem[Sab03]{S}
Joshua~M. Sabloff.
\newblock Invariants of {L}egendrian knots in circle bundles.
\newblock {\em Commun. Contemp. Math.}, 5(4):569--627, 2003.

\end{thebibliography}

\end{document}